\documentclass{amsart}
\usepackage{amsmath, amssymb,amscd}
\usepackage{latexsym, epsfig,color}

\newcommand{\ybG}{y_{\beta(\Gamma)}}
\newcommand{\xbG}{x_{\beta(\Gamma)}}

\newcommand{\C}{\mathbb{C}}
\newcommand{\D}{\mathcal{D}}
\newcommand{\Z}{\mathbb{Z}}
\renewcommand{\P}{\mathbb{P}}
\newcommand{\Q}{\mathbb{Q}}

\newcommand{\GQ}{G_{\Q}}
\newcommand{\Qbar}{\bar{\Q}}
\newcommand{\ra}{\rightarrow}

\newcommand\Aut{\operatorname{Aut}}

\newcommand\Gal{\operatorname{Gal}}

\newtheorem{theorem}{Theorem}

\newtheorem{proposition}[theorem]{Proposition}

\newenvironment{example}{\vspace{2 ex}{\noindent{\bf
Example. }}}{\vspace{2 ex}}

\newenvironment{notationconventions}{\vspace{2 ex}{\noindent{\bf
Notation and Conventions. }}}{\vspace{2 ex}}
\newenvironment{definition}{\vspace{2 ex}{\noindent{\bf
Definition. }}}{\vspace{2 ex}}
\newenvironment{remark}{\vspace{2 ex}{\noindent{\bf
Remark. }}}{\vspace{2 ex}}

\begin{document}

\title[Galois action on dessins d'enfants]{Belyi-extending maps and the Galois action on dessins d'enfants}

\author{Melanie Wood}
          \address{
Department of Mathematics\\
Princeton University\\
Princeton, NJ 08544-1000\\
USA.}
  \date{January 24, 2005}
  \email{ melanie.wood@math.princeton.edu}     
\thanks{This work was supported in part through
NSF grant DMS-9820438, NSA grant MDA904-02-1-0060, and
 VIGRE grant DMS-9983320.}

\begin{abstract}
We study the absolute Galois group by looking for invariants and orbits of its
 faithful action on Grothendieck's dessins d'enfants.  We
define a class of functions called \emph{Belyi-extending maps}, which we use to construct
new Galois invariants of dessins from previously known invariants.  
Belyi-extending maps are the source of the ``new-type'' relations on the injection of the absolute Galois group into the Grothendieck-Teichm\"{u}ller group.
We make explicit how to get from a general Belyi-extending map to formula for its associated invariant which can be implemented in a computer algebra package.  We give an example of a new invariant differing on two dessins which have
the same values for the other readily computable invariants.
\end{abstract}

\maketitle

\section{Introduction}\label{intro}
In this paper we work in the direction of a program to understand $\GQ=\Gal(\Qbar/\Q)$, the absolute Galois group, originally sketched by Grothendieck in
his ambitious research outline \cite{Esquisse}.  One part of the program involves understanding $\GQ$ by studying its action
on certain combinatorial structures, structures so simple at first glance that Grothendieck called them \emph{dessins d'enfants} (children's drawings). 
We construct new invariants of the action of $\GQ$ on dessins d'enfants.  
In fact, we give an example of how an invariant produced from our method distinguishes between two dessin orbits on which previously known invariants agree.  
 
Our basic tools are \emph{Belyi-extending} maps $\beta \colon \P^1 \to \P^1$ that are defined over $\Q$, ramified only over $\{0,1,\infty\}$, and send $\{0,1,\infty\}$ into $\{0,1,\infty\}$.  These maps have the property that when $(X,p)$ is a Belyi pair, $(X,\beta\circ p)$ is also a Belyi pair.  
It follows that each of these $\beta$ acts on the set of dessins.  Moreover, the action of a Belyi-extending map $\beta$ on dessins commutes with the $\GQ$ action on dessins.  This implies that two dessins $\Gamma_1$ and $\Gamma_2$ are in the same Galois orbit when $\beta(\Gamma_1)$ and $\beta(\Gamma_2)$ are.  Thus if $I$ is a $\GQ$-invariant of dessins (e.g. the monodromy group, the rational Nielsen class), so is $I\circ\beta$. 
In \cite{Ell}, Ellenberg defines the class of Belyi-extending maps (without naming them) and notes that 
when $I$ is the monodromy group, the maps give new invariants as described above.  He uses this perspective to construct the cartographic group, but does not go any further with it.

Building $\beta(\Gamma)$ from $\Gamma$ turns out to be quite straightforward.  In Section~\ref{GeomCons} we give a simple, concrete way to construct $\beta(\Gamma)$ geometrically from $\beta$ and $\Gamma$.  Then in Section~\ref{CombCons}, we see how, given a $\beta$, we can easily write down a formula that will produce $\beta(\Gamma)$ from $\Gamma$ for any $\Gamma$ (when the dessins are given as pairs of permutations).  Implementing
these procedures in Maple allowed us to find an example where one of our new invariants distinguishes between two dessins orbits previously indistinguished.  

The action of $\GQ$ on dessins can be refined to an action on the algebraic fundamental group of
$\P^1_\C \setminus \{0,1,\infty\}$.  This allows one to inject $\GQ$ into the Grothendieck-Teichm\"{u}ller group
(see \cite{Ihara2}).  Belyi-extending maps
 can be used to get relations (called ``new-type'' relations) on the image of this injection, analogously to the way
 the new dessin invariants are produced.   
In \cite{Nakamura-Schneps} and \cite{Nakamura-Tsunogai} examples of these relations were produced using specific, very natural Belyi-extending maps.  In \cite{Nakamura}, Nakamura discusses the general properties necessary for a map
to give new-type relations via the Belyi-extending procedure, computes these relations in some specific examples,
and gives a list of examples left for further work.  Belyi-extending maps that are composites of
$PGL(2)$-transforms of the power-raising maps have been
systematically studied in \cite{Anderson-Ihara1}, \cite{Anderson-Ihara2}, and \cite{Ihara3}.

\section{Brief Review of Background and Definitions}\label{Back}

\subsection{Belyi's Theorem and Dessins d'Enfants}\label{Dessindefs}
We start with the following remarkable theorem.

\begin{theorem}[Belyi's Theorem, \cite{Belyi}]
A compact Riemann
surface $C$ is the analytic space associated to an algebraic curve defined
over $\bar\Q$ if and only if there exists a holomorphic map $p_C$ from $C$ to $\P_{\C}^1$ with ramification only over three points.
\end{theorem}

Such a map $p_C$ is called a \emph{Belyi map} and $(C,p_C)$ is called a \emph{Belyi pair}.  We usually (and in this paper) assume without loss of generality that the three points are $\{0,1,\infty\}$.  Two Belyi pairs $(C,p_C)$
 and $(C',p_{C'})$ are isomorphic if there is a Riemann surface isomorphism $f \colon C \to C'$ such that
 $p_{C'} \circ f = p_C$.
Grothendieck found Belyi's Theorem both ``deep and disconcerting'' as it showed a ``profound identity between the combinatorics of finite maps on the one hand, and the geometry of algebraic curves defined over number fields on the other'' \cite{Esquisse}.  He defined dessins d'enfants to capture the combinatorial side of this identity, and we will define them here and see how this identity unfolds.

\begin{definition}
A \emph{dessin d'enfants} (or a \emph{dessin}) is a connected graph with a fixed bipartite structure (vertices labelled 0 or 1 so that each edge has a vertex labelled 0 and a vertex labelled 1) and with fixed cyclic orderings of the edges at each vertex.
\end{definition}

Given a Belyi pair, we can define a dessin on the Riemann surface $C$ as the inverse image via $p_C$ of the real interval [0,1].  The natural
 orientation of $C$ determines a cyclic ordering on the edges incident to each vertex.
 The dessin divides $C$ into faces that correspond to pre-images of $\infty.$  
There is a one-to-one correspondence between dessins and Belyi pairs up to isomorphism.  
A dessin has well-defined faces because of the cyclic orderings at each vertex, and we can glue in these faces to give the Belyi pair, which is unique up to isomorphism by the Riemann Existence Theorem.  
For a dessin $\Gamma$, let the Riemann surface of the associated Belyi pair be called $X_\Gamma$.

\begin{notationconventions}
Throughout this paper, $\P^1$ denotes $\P_{\C}^1$ and $U$ denotes $\P^1 \setminus \{0,1,\infty\}$.  We compose paths in fundamental groups so that $\delta_2\delta_1$ is the path that follows $\delta_1$ and then $\delta_2$. Through Section~\ref{BEInv}, we will let the open real edge $\overline{01}$ be the base-point of the fundamental group of $U$.  To generate $\pi_1(U,\overline{01})$, we choose $x$, a counter-clockwise loop around $0$, and $y$, a counter-clockwise loop around $1$.  We let $F_2$ be the free group on $x$ and $y$, which gives us an isomorphism $\pi_1(U,\overline{01})\approx F_2$.   In $\pi_1(U,\overline{01})$ and $F_2$, we let $z$ be the element such that $xyz=1$.  In $\pi_1(U,\overline{01})$, $z$ is a loop that encircles $\infty$ once.  
\end{notationconventions}  

This notation allows us to describe several other sets that also correspond to dessins and Belyi pairs.  
The following sets are all in one-to-one correspondence:
\begin{enumerate}
\item isomorphism classes of Belyi pairs with Belyi map of degree $d$,
\item isomorphism classes of dessins d'enfants with $d$ edges,
\item isomorphism classes of degree $d$ connected covering spaces of $U$,
\item conjugacy classes of index $d$ subgroups of $\pi_1(U,\overline{01})$,
\item conjugacy classes of homomorphisms of $F_2$ onto a transitive subgroup of the symmetric group $S_d$,
\item conjugacy classes of ordered pairs of permutations that generate a transitive subgroup of $S_d$.
\end{enumerate}
We have seen the correspondence between (1) and (2), and the correspondence between (5) and (6) is clear.  The correspondence between (3) and (4) follows from the standard theory of covering spaces.  To see that (1) and (3) correspond, we consider a Belyi pair $p_C \colon C \to \P^1$.  We can remove $\{0,1,\infty\}$ from $\P^1$ and $p_c^{-1}(\{0,1,\infty\})$ from $C$ to get a covering space $p_c \colon C\setminus p_c^{-1}(\{0,1,\infty\}) \to U$, and we can recover the Belyi pair by compactification of the covering space and $U$.  To see that (2) and (6) correspond, we note from an ordered pair of permutations in $S_d$ we can make a dessin with edges $\{1,2,\dots,d\}$ by taking the cycles of the first permutation to be the cyclic edge orderings around the vertices labelled $0$ and the cycles of the second permutation to be the cyclic edge orderings around the vertices labelled $1$.  This determines a graph with fixed bipartite structure and cyclic edge orderings, and the permutations generate a transitive group exactly when the graph is connected, i.e.\ a dessin.  Often when we refer to a dessin, we implicitly also refer to the corresponding Belyi pair, covering space, class of subgroups, homomorphism of $F_2$, and permutation pair.

\subsection{Action of $\GQ$}\label{GQact}

Dessins d'enfants are particularly interesting because the absolute Galois group acts on them.  Since $\GQ$ acts on elements of $\bar{\Q}$, it acts on the coefficients of polynomials over $\Qbar$, and thus on smooth, algebraic curves defined over $\Qbar$ with maps to $\P^1$ ramified only over $\{0,1,\infty\}$. 
This action is well-defined on dessins.
  Moreover, Grothendieck \cite{Esquisse} noted
that this action is faithful.  It has been since shown by Matzat \cite{Matzat} and Lenstra \cite{Schneps1:1} that $\GQ$ acts faithfully on smaller sets of dessins such as the dessins corresponding to Riemann surfaces of genus $1$ and \emph{dessin trees}, dessins whose underlying graph structure is that of a tree.  Thus, in theory we could understand $\GQ$ by understanding how it acts on dessins, or even dessin trees.  Understanding the $\GQ$ action on dessins is one
step towards understanding $\GQ$.

\subsection{$\GQ$-invariants of Dessins}\label{GQinv}

The next question is how we should study this action, and one plan is to look at the orbits of the action.  A goal of current research is to find invariants of the $\GQ$ action on dessins.  In this paper, we refer to functions on dessins that are constant on $\GQ$-orbits as Galois or $\GQ$-invariants.
The theoretical ideal would be to find a ``complete'' list of invariants, in the sense that any two distinct $\GQ$-orbits could be distinguished by some known invariant.  The search is in particular for invariants that are easily computed from the combinatorial structure of the dessin,
and here we will describe the well-known $\GQ$-invariants that meet this criteria.

The \emph{valency list} of a dessin is an ordered triple of the sets
of valencies of the vertices labelled $0,1,$ and $\infty,$ respectively.
(The $\infty$ valencies are just half the number of edges of each face.)
The valency list of a dessin is a good first invariant.  From it coarser invariants such as the degree of the Belyi map and the genus of the Riemann surface in the Belyi pair can be computed.  The automorphism group of a dessin 
is isomorphic to the automorphism group of the corresponding Belyi pair, and
 is a $\GQ$-invariant.  The automorphisms of the Belyi pair are defined over $\Qbar$ and thus are acted on by $\GQ$ in a way that is compatible with the $\GQ$ action on Belyi pairs.

\subsubsection*{The Monodromy Group}\label{Monodromy}

We can understand the correspondence between dessins and permutation pairs further by considering covering spaces.
A connected covering space $h \colon Z \to U$ of degree $d$, corresponding to a dessin $\Gamma$, determines a transitive homomorphism 
$$\phi_\Gamma \colon F_2 \to S_d$$
up to conjugation in $S_d$.  
This homomorphism is given by the left action of $\pi_1(U,\overline{01})$ on the $d$-element set $h^{-1}(\overline{01})$.  This action is called the \emph{monodromy action} and defined as follows.  Let $e \in h^{-1}(\overline{01})$ and let $\gamma \in \pi_1(U,\overline{01})$.  Then if $\tilde{\gamma}$ is the lift of $\gamma$ to $Z$ with $\tilde{\gamma}(0)$ on the edge $e$, we define $\phi_\Gamma(\gamma)$ to act on $e$ on the left by taking it to the element of $h^{-1}(\overline{01})$ that contains $\tilde{\gamma}(1)$.  Intuitively, we start at $e$ and then follow the lift of the path $\gamma$ to get $\gamma e$.  
It is clear that the action of $x$ and $y$ on $h^{-1}(\overline{01})$, the edges of $\Gamma$,
are the permutations that correspond to the cyclic edge orderings at $0$ and $1$ respectively.  

If $\Gamma$ is a dessin  corresponding to the subgroup $H$ of $F_2,$ then the \emph{monodromy group}
of $\Gamma,$ $M(\Gamma),$ is the image of $\phi_\Gamma$.  Equivalently,  
$$M(\Gamma)= F_2 \left/ \bigcap_{g\in F_2} gHg^{-1}\right. .$$
The monodromy group is also a $\GQ$-invariant of dessins and is easily computed since it is generated by the permutations corresponding to the cyclic edge orderings at $0$ and $1.$

\subsubsection*{The Rational Nielsen Class}

The \emph{Nielsen class} $n(\Gamma)$ of a dessin $\Gamma$ is the data
$$
(M(\Gamma); \phi_\Gamma(x), \phi_\Gamma(y), \phi_\Gamma(z)).
$$
The Nielsen class is defined up to abstract isomorphism of pairs $(G;a,b,c)$ for a group $G$ and elements $a$, $b$, and $c$
that generate $G$ and satisfy $abc=1$. 
Denote the conjugacy class of  $s$ by $[s]$ and let $\widehat{\Z}^\times$ be the invertible
elements in the profinite completion of $\Z$.
The \emph{rational Nielsen class} of a dessin $\Gamma$ is
$$
N(\Gamma)=\{\,(M(\Gamma); [\phi_\Gamma(x)^\lambda], [\phi_\Gamma(y)^\lambda], [\phi_\Gamma(z)^\lambda]) \mid \lambda\in\widehat{\Z}^\times\,\}.
$$
The rational Nielsen classes are also defined up to abstract isomorphism (isomorphisms need not preserve which tuple corresponds to which $\lambda$).
The rational Nielsen class is another $\GQ$-invariant.

\section{Invariants from Belyi-extending Maps}\label{BEInv}

\subsection{New Invariants}\label{NewInv}
To find new invariants we consider a special sort of map as suggested by Ellenberg in \cite[Section 2, Subsection: Cartographic group, and variants]{Ell}.

\begin{definition}
A map $\beta \colon \P^1 \ra \P^1$ is called \emph{Belyi-extending} if $\beta$ is a Belyi map defined over $\Q,$ and 
$\beta(\{0,1,\infty\})\subset\{0,1,\infty\}.$  
\end{definition}

Belyi-extending maps have recently been used to find ``newtype'' equations satisfied by the absolute Galois group
 in terms of its action on the fundamental group of $\P^1\setminus\{0,1,\infty\}$ (see \cite{Nakamura-Tsunogai}, \cite{Nakamura}).  Thus it is natural to investigate whether these maps give new information about the action of
the absolute Galois group on covering spaces of $\P^1\setminus\{0,1,\infty\}$, i.e., dessins.

Note that if $h$ is a Belyi map and $\beta$ is Belyi-extending, then $\beta\circ h$ is also a Belyi map.
If $\beta$ is a Belyi-extending map and $\Gamma$ is a dessin
corresponding to the Belyi pair $(X_\Gamma,p_{X_\Gamma})$, define $\beta(\Gamma)$ to be the dessin corresponding to
$(X_\Gamma,\beta \circ p_{X_\Gamma})$.

\begin{proposition}\label{P:beinv}
Let $(X_\Gamma,p_{X_\Gamma})$ be a Belyi pair corresponding to a dessin $\Gamma$
and $\beta$ be a Belyi-extending map.  
If $I$ is a $\GQ$-invariant of dessins, so is $I \circ \beta$.
\end{proposition}

\begin{proof}
Note that for $\sigma\in\GQ$,
$\sigma(\beta)=\beta$, and thus $\sigma((X_\Gamma,\beta\circ p_{X_\Gamma}))=(\sigma(X_\Gamma),
 \beta\circ\sigma(p_{X_\Gamma})).$  Thus,
$\sigma(\beta(\Gamma))=\beta(\sigma(\Gamma)),$ and so for an invariant $I$, 
$I(\beta(\Gamma))=I(\sigma(\beta(\Gamma)))=I(\beta(\sigma(\Gamma))),$
which proves the proposition.
\end{proof}

As Ellenberg notes in \cite{Ell}, the proof of Proposition~\ref{P:beinv} is the same as the proof he gives in the specific case $\beta(t)=4t(1-t)$ and $I$ being the monodromy group.  Since the monodromy group is a useful and easily computable invariant, we use it to define a whole new class of invariants.

\begin{definition}
If $\beta$ is a Belyi-extending map and $\Gamma$ is a dessin, then we define the $\GQ$-invariant $M_\beta(\Gamma)$ to be $M(\beta(\Gamma))$, the monodromy group of $\beta(\Gamma)$.
\end{definition}

Note that $M =M_{\mathrm{id}}=M_t$.  Also the invariant $M_{4t(1-t)}$ is defined as the \emph{cartographic group}
 in \cite{JonesStreit}.
 Along with these $M_\beta,$ we could consider the automorphism groups, rational Nielsen classes, 
and other invariants of $\beta(\Gamma)$ as new invariants for $\Gamma.$  In this paper we will further investigate the $M_\beta.$

\subsubsection*{Two New Cartographic Groups}

Consider the map $\alpha(t)=4t(1-t)$; $M_\alpha$ is
the cartographic group.  It turns out that $\alpha_1=\alpha\circ(1/t)$ and $\alpha_2=\alpha\circ(t/(t-1))$
lead to new Galois invariants we call the \emph{$\overline{1\infty}$-cartographic group}, $M_{\alpha_1}$, and the
 \emph{$\overline{0\infty}$-cartographic group}, $M_{\alpha_2}$.  These invariants distinguish some orbits indistinguished by the  monodromy group, rational Nielsen class, automorphism group, and traditional cartographic group invariants.  (Take most examples, such as \cite[Ex. 5]{JonesStreit}, of the cartographic group separating orbits and exchange the roles of $0$, $1$, and $\infty$.)  Given the symmetry of the situation, these are natural group invariants to define after defining the cartographic group, but they have
 not previously been mentioned in the literature to the author's knowledge.

\subsection{Geometric Construction of $\beta(\Gamma)$ from $\Gamma$}\label{GeomCons}

We have some examples of easily computable $\GQ$-invariants of dessins and to extend these to useful invariants
via Belyi-extending maps, we will show how to effectively construct $\beta(\Gamma)$ from $\Gamma$.
For a dessin $\Gamma$ drawn as a graph with cyclic edge orderings, we describe how to construct $\beta(\Gamma)$ in a geometric fashion.  We will then see that the geometric construction can be translated
into a combinatorial method.  

For a dessin $\Gamma$, we can form an associated tripartite graph with cyclic edge orderings, which we call
 $T(\Gamma)$.
To construct $T(\Gamma)$, we start with $\Gamma$ and add a vertex labelled $\infty$ for each face of $\Gamma$, 
connecting the vertex with edges to all of the face's vertices in the cyclic order they occur.  (Recall that faces
and cyclic orderings of edges around faces are determined by the cyclic edge orderings at the vertices of
the dessin.)  We can view the real line of $\P^1$ as a triangle with vertices at $0$, $1$, and $\infty$.
Then if $\Gamma$ is associated to the Belyi map $p_{X_\Gamma} \colon X_\Gamma \to \P^1$, $T(\Gamma)$ is simply
the inverse
 image of that triangle via $p_{X_\Gamma}$. 

\begin{figure}[!ht]
\begin{center}
\scalebox{.4}{\includegraphics{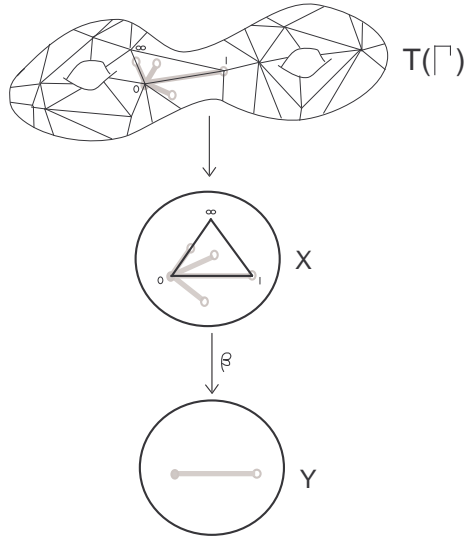}}
\caption{Geometric Construction of a new dessin from a Belyi-extending map}\label{geompic}
\end{center}
\end{figure} 

For a dessin $\Gamma$, we draw $T(\Gamma)$ as a black graph in $X_\Gamma$ over a black $01\infty$
triangle in $X=\P^1$. See Figure~\ref{geompic}.  Now we color gray the dessin in $X$ associated to the Belyi map  $\beta \colon X \to Y=\P^1$, where the $\overline{01}$ edge of $Y$ is also colored gray.  This
gives an overlay of a gray dessin and the black triangle in $X$, which we call the \emph{extending pattern} of $\beta$.  The extending pattern contains the combinatorial information of a super-imposed gray dessin and labelled black triangle.  This information includes the vertices, edges, and cyclic edge orderings. In Figure~\ref{geompic}, we mark the $0$ and $1$ vertices of the gray dessin with filled and open circles, respectively.

We can slice $X$ along the real segments $\overline{0\infty}$ and $\overline{1\infty}$ to get 
a diamond $\D$ composed of the closed upper hemisphere and the open lower hemisphere attached along the edge $\overline{01}$.
\begin{figure}[!ht]
\begin{center}
\scalebox{.5}{\includegraphics{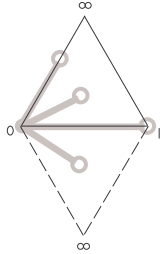}}
\caption{Extending pattern of a dessin}\label{diamond}
\end{center}
\end{figure} 
Figure~\ref{diamond} shows $\D$ in black. Note that $\D$ maps continuously and surjectively onto $X$.  
The pre-images of $\D$ in $X_\Gamma$ partition the surface into diamonds, and the boundaries of the diamonds
are the $\overline{0\infty}$ and $\overline{1\infty}$ edges of $T(\Gamma)$.  Each diamond is composed of a closed triangle over the upper hemisphere of $X$ and an open triangle over the lower hemisphere of $X$, attached at an edge over $\overline{01}$. We can view the extending pattern of $\beta$ in $\D$
as in Figure~\ref{diamond}.

To construct $\beta(\Gamma)$ as an overlay on 
$T(\Gamma)$, we copy the gray dessin from the extending pattern of $\beta$ in $\D$ to each diamond in $X_\Gamma$.  In Figure~\ref{geompic}, we have started this construction by copying the extending pattern of $\beta$ into one diamond of $X_\Gamma$.
  Doing this for every diamond of $X_\Gamma$ will draw $\beta(\Gamma)$ as a gray dessin on top of $T(\Gamma)$ and only depends on the combinatorial information of $T(\Gamma)$ and the extending pattern of $\beta$.   
Thus, $\beta(\Gamma)$ is determined by $\Gamma$ and the extending pattern of $\beta$.  We call two Belyi-extending maps \emph{equivalent} if they have the same extending pattern.

\begin{remark}
If $\beta_1$ and $\beta_2$ have extending patterns that are the same except that the gray labels $0$, $1$, and $\infty$ have been permuted, then $\beta_1=h\circ\beta_2$, where $h$ is a complex
 automorphism of $(\P^1,\{0,1,\infty\})$ (one of $t,1-t,1/t,1/(1-t),(t-1)/t,t/(t-1)$).
It should be noted that for such an $h$, $M=M\circ h$ and thus $M_{\beta_1}=M_{\beta_2}$. 
If $\beta_1$ and $\beta_2$ have extending patterns that are the same except that the labels of the black triangle have been permuted, then $\beta_1=\beta_2\circ h$, where $h\in\Aut(\P^1,\{0,1,\infty\})$. 
This is the case for the traditional cartographic group, our new $\overline{1\infty}$-cartographic group, and our new
$\overline{0\infty}$-cartographic group.
Since often the Belyi-extending $\beta$ will not be symmetric with respect to
 $\{0,1,\infty\}$, $\beta$ and $\beta\circ h$ can give different information as Belyi-extending maps.  For example, when $\beta_1=\beta_2\circ h$,
 the invariants $M_{\beta_1}$ and $M_{\beta_2}$ usually give different 
 information about $\GQ$-orbits of dessins, as with the cartographic groups.  
    Heuristically, if the invariant $M_\beta$
was somehow sensitive to the 0 vertices of a dessin, $M_{\beta\circ(1-t)}$ would be sensitive in that way to the 1 vertices.  
\end{remark}

\subsection{Combinatorial Construction of $\beta(\Gamma)$ from $\Gamma$}\label{CombCons}

Combinatorially, we have the extending pattern of $\beta$ and we want a formula in which we can plug the monodromy
permutations of $\Gamma$ to get the monodromy permutations of $\beta(\Gamma)$.  
Recall that
$$\phi_\Delta \colon \pi_1(\P^1 \setminus \{0,1,\infty\},\overline{01})\approx F_2 \to S_d$$
gives the monodromy action on the edges of the dessin $\Delta$.
In this section, let
$\phi_\Delta(w)=w_\Delta$.  The extending pattern of $\beta$ gives us the permutation pair $[x_\beta,y_\beta]$
corresponding to the Belyi pair $(\P^1,\beta)$, and let $E_\beta$ be the set of edges of the corresponding dessin.
 Let $E_\Gamma$ be the set of edges of a dessin $\Gamma$, which is just
 $\{1,2,\dots,\deg \Gamma\}$ (since we will not use information about $\Gamma$ to construct our formula for
 $\beta(\Gamma)$). 

\begin{proposition}\label{labelascross}
Given the geometric construction of Section~\ref{GeomCons}, we can naturally identify the
edges of $\beta(\Gamma)$ with the elements of $E_\Gamma \times E_\beta$.
\end{proposition}

\begin{proof}
From Section~\ref{GeomCons}, $X_\Gamma$ can be partitioned into diamonds, and
each diamond
can be naturally identified with the element of $E_\Gamma$ that it contains.
Also, each of these diamonds contains exactly one pre-image of every element of $E_\beta$.  Thus, each edge of $\beta(\Gamma)$ 
can uniquely be identified by the diamond it is in and the element of $E_\beta$ it is over.
\end{proof}

First we will find the permutation $\xbG$, and then note $\ybG$ can be found analogously.  
Let $(a,b)\in E_\Gamma \times E_\beta$ be an edge of $\beta(\Gamma)$, and we will
find $\xbG(a,b)$.  The monodromy elements $\xbG$ and $x_\beta$ come from the same
loop $x$ in the fundamental group
of the downstairs $\P^1\setminus\{0,1,\infty\}$.  Thus the action of $\xbG$ in the covering space
associated to $\beta(\Gamma)$ restricts to the action of $x_\beta$ in the covering space
associated to $(\P^1,\beta)$.  This means that $\xbG$ acts as $x_\beta$ on the $E_\beta$ factor of the edges of
$\beta(\Gamma)$.

Now consider the $0$ vertex of the edge $b$ in the extending pattern of $\beta$ (the vertex of $b$ that is above $0$, not the $0$ vertex of the black triangle).  Note if
in the cyclic edge orderings around this $0$ vertex, the black $\overline{0\infty}$ segment or
the black $\overline{1\infty}$ segment is between the edges $b$ and $x_\beta b$.  We wish
to find which of these black segments, if any, the monodromy path $x_\beta$ that takes $b$ to $x_\beta b$
crosses.  We say the path \emph{crosses} one of these segments if it goes from the closed upper hemisphere of $X$
to the open lower hemisphere of $X$ or vice versa by way of intersecting the segment.
From this information, we can find $w(b)$, the action of $\xbG$ on the $E_\Gamma$ factor, so that $\xbG(a,b)=(w(b)a,x_\beta b)$.
If, in going counter-clockwise around the 0 vertex of $b$ to get from $b$ to $x_\beta b$, the monodromy path crosses

\begin{enumerate}
\item the black $\overline{0\infty}$ from the upper to lower hemisphere of $X$, then
$w(b)=x_\Gamma$,
\item the black $\overline{0\infty}$  from the lower to upper hemisphere of $X$, then
$w(b)=x_\Gamma^{-1}$,
\item the black $\overline{1\infty}$ from the lower to upper hemisphere of $X$, then
$w(b)=y_\Gamma$,
\item the black $\overline{1\infty}$ from the upper to lower hemisphere of $X$, then
$w(b)=y_\Gamma^{-1}$,
\item neither the black $\overline{0\infty}$  nor
the black $\overline{1\infty}$, then $w(b)$ is the identity, and
\item two of these edges, then $w(b)$ is the product of the elements of $M(\Gamma)$ given above for those edges, such that the element corresponding to the first edge crossed is on the right.
\end{enumerate}

Often, and always when the relevant vertex is in the interior of the diamond, the action will
fall into case (5).  The results of the cases given are clear because when a path crosses
a black $\overline{0\infty}$ or $\overline{1\infty}$, it goes from one diamond in $X_\Gamma$ to
another and the relationship between the first diamond and the second is exactly as
listed above. So given the extending pattern of $\beta$, we can determine exactly which elements of
$E_\beta$ fall into each of the six cases above.  This allows us to write down a formula
for $\xbG$ in terms of $x_\Gamma$, $y_\Gamma$, and $x_\beta$.  Similarly we can write
 $\ybG$ in terms of $x_\Gamma$, $y_\Gamma$, and $y_\beta$.
 
\begin{remark}
There is a unified story to the multiple cases given above.  In each case, the action
of $\xbG$ on the $E_\Gamma$ factor is an element $w(b)$ that depends on $b$,
the value in the $E_\beta$ factor.  To determine $w(b)$, consider the path $\rho(b)$ in $X=\P^1$ that starts
at the edge $b$ and goes counterclockwise around the $0$ vertex of $b$ very near the vertex and ends
at $x_\beta b$.  The path $\rho(b)$ can be viewed as an element of $\pi_1(X\setminus\{0,1,\infty\},\D)$,
where the base-point is the map of $\D$ into $X$ described in Section~\ref{GeomCons}.  However,
we can canonically identify $\pi_1(X\setminus\{0,1,\infty\},\D)$ with 
$\pi_1(X\setminus\{0,1,\infty\},\overline{01})$ since $\overline{01}$ lies inside $\D$.  This allows
us to identify the path $\rho(b)$ with an element $W(b)\in\pi_1(X\setminus\{0,1,\infty\},\overline{01})$.  
Then $w(b)$ is just the image of $W(b)$ in $M(\Gamma)$.
\end{remark}

Some examples will illustrate the combinatorial construction.

\begin{example}\label{cf}
Consider the map $\gamma(t)=-27(t^3-t^2)/4.$  Note that $\gamma(\{0,1,\infty\})=\{0,\infty\},$ $\gamma$ is ramified only over
$\{0,1,\infty\}$, and $\gamma$ is defined over $\Q$.  Thus $\gamma$ is Belyi-extending. 
\begin{figure}[!ht]
\begin{center}
\scalebox{.5}{\includegraphics{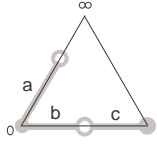}}
\caption{Extending pattern of $\gamma$}\label{gammapat}
\end{center}
\end{figure}
 We can draw 
the extending pattern of $\gamma$, as shown in Figure~\ref{gammapat}.
We see that $[x_\gamma,y_\gamma]=[(a~b), (b~c)]$.
Let the edges of $\Gamma$ be denoted $1,2,\dots,d$.  Then $x_{\gamma(\Gamma)}$ is composed of all the cycles of the form
\begin{eqnarray*}
&\left( (k,a)~(x_\Gamma k,b)~(x_\Gamma k,a)~(x_\Gamma^2 k,b)~(x_\Gamma^2 k,a)\cdots \right)\\
\mbox{and } &\left( (k,c)~(y_\Gamma k,c)~(y_\Gamma^2 k,c)\cdots \right)
\end{eqnarray*}
for $1\leq k\leq d$, and $y_{\gamma(\Gamma)}$ is composed of all the cycles of the form
$$
\left( (k,b)~(k,c) \right),
$$
for $1\leq k\leq d$.

Thus we can write a simple procedure that will, for the permutation pair of any $\Gamma$, produce the permutation
pair of $\beta(\Gamma)$.  With this we can make computations in a program such as Maple about $M_\gamma(\Gamma)$ from the permutation pair
of $\Gamma.$  We can also easily compute $M_{\gamma\circ h}$ for any  $h\in\Aut(\P^1,\{0,1,\infty\})$
by exchanging the roles of $x_\Gamma$, $y_\Gamma$, and $(x_\Gamma y_\Gamma)^{-1}$ in the procedure.
\end{example}

\begin{example}\label{S3example}
Consider the map $\xi(t)=27t^2(t-1)^2/(4(t^2-t+1)^3),$ and note that it is Belyi-extending.  This is the quotient map by the automorphisms of $(\P^1,\{0,1,\infty\})$.  

\begin{figure}[!ht]
\begin{center}
\scalebox{.5}{\includegraphics{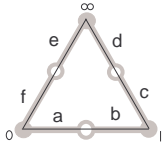}}
\caption{Extending pattern of $\xi$}\label{xipat}
\end{center}
\end{figure}

The extending pattern of $\xi$ is shown in Figure~\ref{xipat}. 
We see that $[x_\xi,y_\xi]=[(a~f)(b~c)(d~e), (a~b)(c~d)(e~f)]$.
Let the edges of $\Gamma$ be denoted $1,2,\dots,d$ and let $z_\Gamma=(x_\Gamma y_\Gamma)^{-1}$.  Then $x_{\xi(\Gamma)}$ is composed of all the cycles of the forms
$$
\left( (k,a)~(k,f)~(x_\Gamma k,a)~(x_\Gamma k,f)~(x_\Gamma^2 k,a)~(x_\Gamma^2 k,f)\cdots \right),$$
$$
\left( (k,c)~(k,b)~(y_\Gamma k,c)~(y_\Gamma k,b)~(y_\Gamma^2 k,c)~(y_\Gamma^2 k,b)\cdots \right), \mbox{ and}
$$
$$
\left( (k,e)~(k,d)~(z_\Gamma k,e)~(z_\Gamma k,d)~(z_\Gamma^2 k,e)~(z_\Gamma^2 k,d)\cdots \right),
$$
for $1\leq k\leq d$, and $y_{\xi(\Gamma)}$ is composed of all the cycles of the forms
$$
\left( (k,a)~(k,b) \right),
\left( (k,c)~(k,d) \right), \mbox{ and}
\left( (k,e)~(k,f) \right),
$$
for $1\leq k\leq d$.
 
Note that since $\xi=\xi\circ h$ for any $h\in\Aut(\P^1,\{0,1,\infty\})$,we get no new invariants 
from considering $\xi\circ h$. 
\end{example}

The new invariants that we have produced with Belyi-extending maps
can give information about the $\GQ$-orbits of dessins that we can't get from old invariants.

\begin{proposition}\label{xigood}
Let $\Delta$ be the dessin corresponding to the permutation pair
$$[(1~2~3~4)(5~6~7)(8~9), (1~8~4~7)(2~3~10)(5~6)]$$
and $\Omega$ be the dessin corresponding to the  permutation pair
$$[(1~2~3~4)(5~6~7)(8~9), (1~3~8~9)(2~10)(4~5~6)]$$
(shown in Figure~\ref{mydes}).
\begin{figure}
\begin{center}
\begin{picture}(0,0)%
\includegraphics{delta-omega.pstex}%
\end{picture}%
\setlength{\unitlength}{1480sp}%
\begingroup\makeatletter\ifx\SetFigFont\undefined%
\gdef\SetFigFont#1#2#3#4#5{%
  \reset@font\fontsize{#1}{#2pt}%
  \fontfamily{#3}\fontseries{#4}\fontshape{#5}%
  \selectfont}%
\fi\endgroup%
\begin{picture}(7366,3924)(818,-3223)
\put(7501, 89){\makebox(0,0)[lb]{\smash{{\SetFigFont{8}{9.6}{\rmdefault}{\mddefault}{\updefault}{\color[rgb]{0,0,0}$\Delta$}%
}}}}
\put(1126,-2611){\makebox(0,0)[lb]{\smash{{\SetFigFont{8}{9.6}{\rmdefault}{\mddefault}{\updefault}{\color[rgb]{0,0,0}$\Omega$}%
}}}}
\end{picture}%
\caption{}\label{mydes}
\end{center}
\end{figure}
Then $\Delta$ and $\Omega$ have the same valency lists, automorphism groups,
monodromy groups, cartographic groups, and rational Nielsen classes.  However,
the $M_\xi$ groups of $\Delta$ and $\Omega$ differ in size.  Thus, $\Delta$ and $\Omega$
are in different $\GQ$-orbits
\end{proposition}

\begin{proof}
All of the necessary computations can be done with Maple.  We can compute that
$\phi_\Delta(z)=(1~6~7~3)(2~10)(4~9~8)$ and $\phi_\Omega(z)=(1~8~2~10)(6~7)(3~4~5)$, and
thus see that $\Delta$ and $\Omega$ have the same valency lists.  We can compute with Maple,
or notice by drawing the dessins, that both have trivial automorphism group. We can notice this
easily since both have a single 0 vertex of degree 1, which means that any automorphism fixes the incident edge and
thus all edges.
We can compute with Maple the monodromy and cartographic groups of $\Delta$ and $\Omega$ as permutation groups and 
can verify that these invariants are the same for the two dessins.
In fact, just computing the size of the monodromy groups of $\Delta$ and $\Omega$ shows
that each has order $10!/2$.  This implies that both monodromy groups are the alternating
group $A_{10}$ since they are index two subgroups of $S_{10}.$  Notice that in the monodromy group
of $\Delta$ or $\Omega,$ $[g]$ includes all elements of the same cycle type as $g$ for which
the relabelling permutation from that element to $g$ is even. Then note that 
$\phi_\Delta(x)=\phi_\Omega(x)$,
conjugation by $(8~ 3~ 5~ 2~ 4)(7~ 9)(10~ 6)$
relabels $\phi_\Delta(y)$ to $\phi_\Omega(y)$,
and conjugation by $(6~ 8~ 5~ 9~ 4~ 3~ 10~ 7~ 2)$
relabels $\phi_\Delta(z)$ to $\phi_\Omega(z)$.  Thus $\Delta$ and $\Omega$ have the same rational Nielsen class.
However, computation in Maple using the description of $\xi(\Gamma)$ above shows that
$M_\xi(\Delta)$ has size $19752284160000$ while $M_\xi(\Omega)$ has size $214066877211724763979841536000000000000$.
(This example was constructed so that $\Delta$ would have the $S_3$ symmetry that $\xi$ quotients by, but
there is no reason to think this would be the only example where $M_\xi$ provides new information,
 only the easiest to construct.)
\end{proof}

Using the procedure we have described for getting $\beta(\Gamma)$ from $\Gamma$, we can compute that $\Delta$ and $\Omega$ agree on some other $\GQ$-invariants we have defined in this paper.
They agree on the $\overline{0\infty}$-cartographic group, the $\overline{1\infty}$-cartographic group, and 
the groups $M_{\gamma\circ h}$ for any  $h\in\Aut(\P^1,\{0,1,\infty\})$.  There are other known Galois invariants of dessins which we have not checked here, such as Ellenberg's braid group invariant \cite{Ell} and Fried's lifting invariants \cite{Fried}.  However, these invariants are not readily computable from the combinatorial data of the dessin.  Our new invariants from Belyi-extending maps have the benefit of being easily computable with a computer algebra package.

\section{Directions for Further Research}\label{further}

\begin{figure}
\begin{center}
\scalebox{.3}{\includegraphics{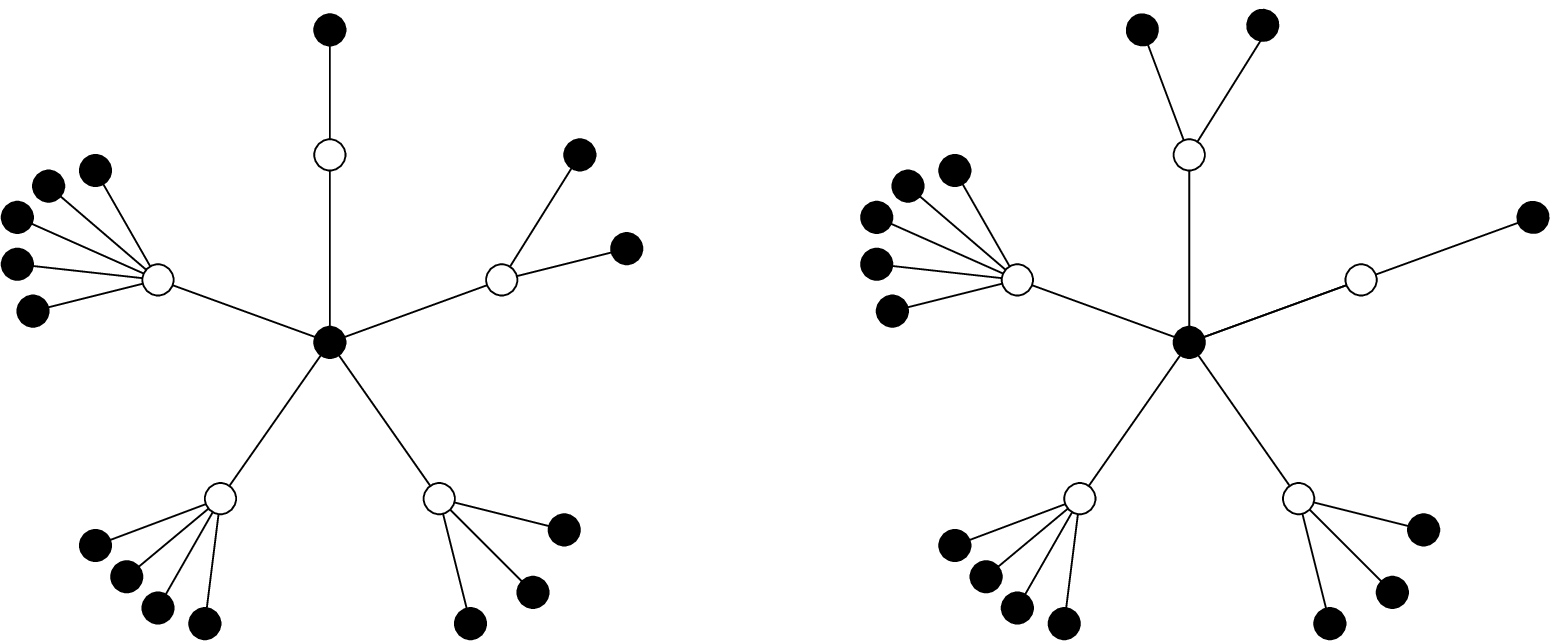}}
\caption{}\label{flowers}
\end{center}
\end{figure}

The most well-known example (\cite[Ex. 1]{Schneps1:1}) of two dessins which have been computed to be in different orbits (by finding
Belyi maps that give the dessins) but are indistinguished by the known Galois invariants are shown in Figure~\ref{flowers}.

However, finding Belyi maps that produce a certain dessin is very difficult computationally, and so usually
when we know two dessins are in different orbits it is because we know an invariant on which they differ.
It is likely there are many dessin orbits which could be distinguished by computing new invariants but which
are not currently distinguished by any means.

    Since Belyi-extending maps give us new, computable $\GQ$-invariants of dessins, a natural direction of
    research is to find such maps.  Recall that we can take such a map and compose on the left or right
    an automorphism of $(\P^1,\{0,1,\infty\})$ to get another such map.  We refer to a set of Belyi-extending maps that arise from each other by composing with automorphisms a \emph{family of Belyi-extending maps}. Considering other maps in the family of $\alpha$ was what led us to the $\overline{0\infty}$- and $\overline{1\infty}$-cartographic groups.  Finding Belyi-extending maps amounts
    to looking for rational solutions to polynomial equations and though this is difficult to do in general,
    it can be done in many specific cases.  Recall that for a Belyi-extending map $\beta$ and a dessin $\Gamma$, the only structure that determines how to get $\beta(\Gamma)$ from $\Gamma$ is the extending pattern of $\beta$ defined in Section~\ref{GeomCons} and that we call two Belyi-extending maps equivalent if they have the same extending pattern.   
    
\begin{figure}[!ht]
\begin{center}
\scalebox{.3}{\includegraphics{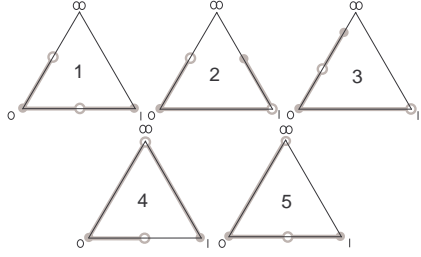}}
\caption{Extending patterns of several Belyi-extending maps}\label{betaspat}
\end{center}
\end{figure}

    The following degree 3 maps are Belyi-extending, non-equivalent, and all in distinct families.  Their extending patterns are shown in Figure~\ref{betaspat}, with the extending pattern of $\beta_i$ labelled $i$.
      \begin{itemize}
    \item $\beta_1(t)=\gamma(t)=-\frac{27}{4}(t^3-t^2)$
    
    \item $\beta_2(t)=-2t^3+3t^2$
    
    \item $\beta_3(t)=\frac{t^3+3t^2}{4}$
    
    \item $\beta_4(t)=\frac{27t^2(t-1)}{(3t-1)^3}$
    
    \item  $\beta_5(t)=t^2(t-1)/\left(t-\frac{4}{3}\right)^3$
   
    \end{itemize}

  The following maps are quotient maps of $\P^1$ by some of its finite automorphism group actions.  The equations for the 
  maps  $\beta_{T}$ and  $\beta_O$ were found in \cite{CouvGran}.  They are all Belyi-extending, non-equivalent, and in distinct families.  It is important to note that there could be other Belyi-extending maps that are quotients by group actions, because it is not just the quotient that determines the equivalence class of the maps but also the position of the real triangle in the domain, i.e. the precise extending pattern and not just the dessin $(\P^1, \beta)$.
 
    \begin{description}
    \item[$\mathbf{C_n}$] $\mu_n(t)=t^n$
    
    \item[$\mathbf{D_n}$] $\nu_n(t)=(t^n+1)^2/(t^n-1)^2$
    
    \item[$\mathbf{D_3=S_3}$] $\xi(t)=\frac{27t^2(t-1)^2}{4(t^2-t+1)^3}$ 
   
    \item[Tetrahedron] $\beta_{T}(t)=\frac{t^3(t^3+8)^3}{(t^6-20t^3-8)^2}$
    
    \item[Octahedron] $\beta_{O}(t)=\frac{108t^4(t^4-1)^4}{(t^8+14t^4+1)^3}$ 
        
    \end{description}
    
The above Belyi-extending maps or other maps in their families could potentially lead to new $\GQ$-invariants of dessins. 
The family of Belyi-extending maps is closed under composition, so composing the above maps gives
infinitely many Belyi-extending maps.
For all of these Belyi-extending maps $\beta$ there are many new invariants to consider.  
We could compose
any $\GQ$-invariant with the action of $\beta$ on dessins to get a new dessin invariant; however, 
the monodromy group, rational Nielsen class,  and
automorphism group give the most easily computable new invariants.  It can be shown that the invariants $M_{\mu_n}$ and the rational Nielsen classes of $\mu_n(\Gamma)$ can be constructed from previously known invariants of $\Gamma$, but we saw that $M_\xi$ is a new invariant and it is not known whether any of the other maps give new invariants.  Investigating this question for various $\beta$ would be a interesting direction of further research.

\section{Acknowledgements}
The author wishes to thank Richard Hain for suggesting that she study the exciting field of dessins d'enfants and advising her along the way, Jordan Ellenberg for helpful conversations, and Makoto Matsumoto for his support, encouragement, and careful, detailed answers to many questions.

\end{document}